\documentclass[12pt]{article}

\usepackage{latexsym,amsmath,amssymb}

\begin{document}

{\Large

\begin{center}

\vskip 0.5cm {\bf \MakeUppercase{A little about models} }

\vskip 0.5cm

{\bf I.V.~Konnov}

\vskip 0.5cm

{\em Kazan, e-mail: konn-igor@ya.ru}
\end{center}

}

{\large

\vskip 0.9cm

 \hfill{\sl Sapienti sat}


\section{The fallacy of simple concepts}

In many sciences, it is customary to create models. Although there are some
sciences where most knowledge is obtained through direct
observations (and previously these sciences prevailed), but
a fundamentally new level can only be achieved by creating models that
clarify the essence of studied processes and phenomena.

Since the author's experience is limited to mathematical models,
 first of all we will talk about them. The role of these models in other sciences, and in other areas
of activity is very significant. For example, according to Leonardo da
Vinci, \lq\lq no human investigation can be said true science, if it cannot be demonstrated mathematically".
However, there are also opposite opinions about
negative influence of mathematical models in other sciences
(see, e.g., \cite{Ali80}). These notes can be regarded as
an attempt to clarify some important aspects of this issue.

What, in fact, is a mathematical model of some real
or imaginary system? A set of relationships that link different
parameters and variables  and reflect all significant relationships among
elements (and subsystems) of this system, as well as the relationships with the system environment.
It is obvious that one always try to create any model
as simple as possible since it will be used to
solve various behavioral assessment tasks related to
the system, to estimate the influence of certain parameters on it, etc.
I would also like to emphasize that the solutions to all these problems are actually
describe properties of the model, although they are usually transferred to
the actual system that the model was built on, and this difference
 is very important for right understanding of the mathematical modeling process.

A special feature of mathematical models is that they are based on
a given formal structure, i.e., enough
clearly defined elements and parts of the system under study and all
possible interrelations among them, so that they can be
written in the form of mathematical relations.

For example, {\em number}  is one of the most fundamental concepts in Mathematics. The use of numbers (count)
or the corresponding scalar variable (or parameter) in essence
means that all properties of the enumerated objects are not important
for the corresponding mathematical model, except for their
quantity. That is, they in fact belong to the same type (class) in this model.
Of course, some other objects that are more complex than numbers can be also used.
 But introducing some common preference relation for these general objects
 also postulates implicitly that they belong to the same class because they allow comparison
based on the entered relation.

Based on this, it is quite clear that a formal structure is being introduced
automatically explicitly or implicitly already when creating any
mathematical model. Of course, such a model can be investigated
independently as a mathematical problem, in particular, to find out
questions about the existence of solutions, their type, and their dependence on parameters,
calculation capabilities, etc. But there is always a question about
the application areas of the model. In other words, one has to define
 what real systems have a structure that
 it is adequately reflected by the structure of the constructed model.
Examples of models that have applications in completely different
areas are fairly well known, but the main question is in
methods for checking the adequacy of models.

Indeed, incorrect structuring (formalization) of the original
systems when creating a mathematical model can
 lead to useless multi-cost work and, moreover,
 to false conclusions. At the same time, the successful application of the model in one area
does not guarantee its adequacy in another area, or
in the case of changing the original real system. In this connection
the situation with the use of models in
Physics, where mathematical model construction has been carried out for centuries,
is rather favorable. Indeed, in many cases it was not just possible to create suitable models of
processes and phenomena, but also specify the conditions (ranges) of their successful
applications (see, e.g., \cite{SM05,KP83}). In other sciences, the situation is
not so favorable in many respects due to the presence of poorly formalized
systems (see \cite{EM79}). It should be added that the model that is not structured correctly
cannot be improved by using objects with more general properties or more complicated techniques.

A well-known example of incorrect formalization is
 Ptolemy's geocentric system.  Its application to position determination of
stars and planets met constantly detected deviations that forced one to
use more and more corrections in formulas. Noteworthy,
the initial use of the heliocentric system by Copernicus also revealed deviations from the instant
positions, but they were caused by an inaccurate shape adopted for orbits
(circles instead of ellipses).

Note that various general system structures and their properties
are studied separately in the theory of systems, and different approaches are used for this purpose
(see, e.g., \cite{Moi81,VVD83}). However, the issues of compliance of
structures of the model and real system under study, as a rule,
you have to resolve for yourself.


\section{Impact of uncertainty}

Thus, it turns out that the main problems when building
adequate models arise because of the great complexity of the source systems
and the presence of a variety of uncertain factors, which makes it naturally difficult
to determine the appropriate structure of
models. For example, because of the fundamental impossibility to describe
behavior of each separate gas molecule in statistical physics,
one was able to build an adequate model only at the macro level, i.e. to describe, say,
the behavior of a certain volume of the gas as a whole.

Let's illustrate the effect of uncertainty using the {\em
consumer demand model} in the neoclassical theory of market equilibrium
(see, for example, \cite{Eke79}). In this model, the consumer's demand of goods
defined as a set of solutions to the consumer utility maximization problem
on the budget set where prices of goods
are specified as external parameters. The key element of this model
is just this utility function (or, if necessary, the consumer preference relation),
 whose existence is deduced from the consumer's ability to accurately compare values
of all goods. Such a complete determinism of tastes is quite
corresponds to the real behavior of a buyer who went to a local market (or fair) in
the Middle Ages  held, say, every week, where the same
limited set of home-made products with similar properties was proposed for sale
every time and almost the same amount of money was used
for purchases  there  (see \cite{Kon16a}). It is clear that
a major change in prices could be only invoked by certain external influence
 on these conditions. However, the specified \lq\lq marginal"  demand model
used in neoclassical theory to describe general market behavior within
a sufficiently large time interval, with a large
number of participants and a wide variety of products, so that the same
product produced at a different time or place is considered as a different product.
This approach leads to infinite dimensional models of market equilibrium,
which are very difficult even to determine the conditions of existence of
solutions (see, for example, \cite{ABB90}).
It is obvious that
the consumer is not able to accurately evaluate the usefulness of all products
under these conditions.
Aware of this drawback, but trying to keep this model of
consumer behavior, as the main properties of the
perfect competition model are based on it, supporters of the
neoclassical theory proposed to generalize the concept of the
utility function based on the so-called \lq\lq rational expectations", which allowed
one to maintain the basic model structure without changes. Meanwhile, the main
difference with the usual utility function is that the new
does contain undefined factors, and the level of this
uncertainty may be arbitrarily high. Therefore,
the assumption of maintaining the \lq\lq marginal" consumer behavior
in case of unreliable data seems absolutely unrealistic.
For this reason, a fundamentally different type of the model is required
to describe consumer's behavior in general.

A fairly popular approach to research and solving problems with
uncertainty is the use of random values, so that
the uncertainty is simply identified with the
randomness in many works. Recall that any random variable is determined on
a set (or space) of elementary events together with its
probability distribution (normalized measure) on this set,
and may be of a continuous or discrete type. However, setting such
a measure itself and associated values does not mean that
just a random variable is determined, as well as any vector with non-negative coordinates,
the sum of which is equal to one is not automatically a
probability distribution.  This property requires special
conditions that justify the possibility of utilization of a
probabilistic (stochastic) model.

Namely, setting the probability distribution requires {\em
statistical stability}, i.e. evidence based on
multiple observations under the same conditions. Obviously, the
statistical stability is achieved only when observing
a sufficiently homogeneous and independent process or
 phenomenon. On the other hand, the presence of a common measure
 in the form of the probability distribution for elementary
events also indicates that these events are of the same type,
as noted above about the implicit properties of applying numbers (scalar
variables). Thus, the streamlined utilization of probabilistic
models for essentially heterogeneous diverse phenomena is
incorrect.  In addition to the specified conditions, this requires
definition of a formal structure that would be adequate to that of
the source real system. From this we can conclude that not every
uncertainty can be represented by a random variable, and
that these concepts are not equivalent.

A fairly standard technique in game theory describing models of
conflict situations, i.e. models with uncertain factors,
is the utilization of mixed strategies that are nothing but
probability distributions on the set of (ordinary) pure strategies,
which are then elementary events. This approach leads to
a significant complication of the original model, but relaxes sufficient conditions
for existence of equilibrium states (see, for example, \cite{Owe95}).
The presence of equilibrium states makes the behavior of the described
conflict system quite predictable. The implementation of mixed strategies by players
consists in conducting a random experiment in accordance with the corresponding
distribution  and selecting the pure strategy obtained for further actions.
This approach looks quite artificial, moreover, it makes sense only in case
of multiple repetition of playing this game. Then the players' utilities will tend
to their average values, that is, to the game value in mixed strategies.
If the game is played once, or a few times, the usefulness of mixed strategies
becomes doubtful. In some cases
it is possible to utilize the so-called \lq\lq physical mixture" of strategies (see \cite{Ven72}).
For example, if a pure strategy is to select a crop for
sowing a field, there is no need to produce a random
experiment to select it. One can then just sow the field in proportions
specified by the probability distribution. In the general case, players are more likely to
will base their actions on any additional information about the other
 participants, i.e. the game will transform into some multi-stage
 procedure. Detailed discussion of  utilization of mixed strategies
 for finding solutions of various games can be found in \cite[\S 14,16]{Ger71}, \cite[\S 11]{Ger76}.


\section{Information flows}

 The model structure for various socio-economic systems,
 industrial, transport and communication, and in general for
systems related to human activity should include
 information exchange schemes among elements
 and blocks (subsystems).  This is one of the main differences from
 the structure of models in many natural sciences, such as Physics. Under {\em
information} hereafter will only be understood as its
content, rather than its volume recorded on any media.
That is, description of elements and subsystems themselves and their relationships
is not sufficient for adequate definition of the required system structure
without description of the information exchange scheme.

For example, classical models of perfect and imperfect competition
describe two different types of decentralized systems in
Economics. Separate actions of economic agents (elements of the system)
cannot affect the state of the whole system in
the Walrasian type perfect competition models.
But common actions of economic agents can change the system state.
Therefore, each of them in principle need not use information about actions (prices, assortment, volumes) or
interests of some other separate agent.
 Instead, the agents use information about integral indicators of the entire system
 (for instance, good prices), which
may be available to them, although the mechanism for determining common market prices
is not clearly defined in the available models (see, for example,
\cite{Kon15e}). On the other hand,  separate actions of each economic agent in imperfect
competition models can change the state of the entire system
and, in particular, affect any other economic agent, so the participants will use information about
actions and interests of others when choosing their own actions. As a result,
it turns out a fundamentally different game-theoretic model, with a different
information exchange scheme.

One can also find a lot of examples in history when
states with the same government institutions were managed in completely different manners,
 that is, with  different information exchange schemes among these institutions.

In addition to the above (custom) definition,
it is also common to define the information concept
 based on probabilistic representations, which for the difference
 will be denoted as {\em information (p)}.
For example, let us suppose that the set of all possible events for
an investigated system (object) is finite $(n)$ for simplicity, so as
the whole set is the union of these $ n $ elementary events.
A state $S$ of the  system is then defined by using some
probability distribution $p= (p_{1},\ldots,p_{n})$ on this set, that
allows one to calculate its \lq\lq entropy"
$$
H(S)=- \sum \limits_{i=1}^{n} p_{i}\log p_{i},
$$
which is considered as an uncertainty measure of this state. Therefore,
the maximal entropy corresponds to the greatest
uncertainty, i.e. to the state where all the elementary events may occur with the
same probability, or simply where $p_{i}=1/n$, $i=1, \dots, n$. The probability of events,
as usual,  can be conditional or unconditional, and the difference
 between two states is just determined
as the information (p); i.e., decrease of entropy gives the positive
value of the information (p) in this transition (see, for example, \cite{Str75}).
However, this raises a natural question about the generality of such
approach, since it in fact states the possibility of existence of one-dimensional
representations for any diverse processes and phenomena. For example,
reception of some new knowledge about the system invokes transition to a new state,
but this knowledge can be somewhat incorrect, so how does one evaluate
the change of the entropy in this case?
Is it possible to measure the amount of information in any object at all?
It is clear that this measurement is only possible within
some suitable model with the specified formal structure. In general,
 both the objects themselves and events with them contain in fact infinite
amounts of information, so such \lq\lq ubiquitous" measurements are meaningless. For instance,
 when replacing one our car with another of the same kind we clearly understand that
 it is not the same car. But within the framework of our model (representation) of a car, it will perform
 the same necessary functions, and all the differences
 between them, though infinite, are insignificant to us. Hence,
we can simply ignore these differences and consider the new car as the same. That's why
the amount of necessary information about it can be considered finite.
Thus, the information (p) must be related to some
 specific structure of the model used. Note that different model
structures can in principle use the same set of
information about one object.

This concept is also associated with the reflected information (data), i.e.
recorded on certain media, which can be denoted as
{\em information (d)}. There are many tasks associated with
efficient processing, storage, and transmission of the information (d) on
various devices that are not directly related to the
model creation. We can only note that the information (d) appears initially on
the devices within some specific models,
and the processing issues are to some extent related to the models and their
structures in which they will be used. The differences between these
concepts are clearly indicated, for example, in \cite[p. 111]{YY73}:
\lq\lq The concept of information arose directly from the problems of
communication theory and was specially selected to meet the objectives of this theory.
 Since the transmission of a fixed length message over a communication line requires
approximately  the same time and expenses both in the case of an insignificant or
even false message and in the case of a message about
the greatest discovery, we must assume  from the point of view of the  communication
theory that the amount of information in both these messages is also the same".


\section{Additional examples}

For more clarity, we will give additional illustrative examples of models
of quite complex systems.


\subsection{Mean field games}

This model is intended to describe the behavior of a team involving
a sufficiently large number of dynamic active elements (players), i.e.
each of them has its own goal
 function and state equation that determines the relationship
between the trajectories and control functions,
as in the custom differential game. The model is based on
the assumption that the players differ only in random terms.
Then it is suggested to go to the averaged values and after taking
the limit on the number of players one can get an
 optimal control problem, which will
 consist in maximizing the corresponding averaged functional on
 the averaged equation of state. Obviously, this approach
 is based on the direct transfer of modeling principles from
statistical physics. But in this case its utilization for
 socio-economic applications is emphasized (see, for example,
\cite{GLL11}). However, this raises a natural question about the
validity of the approach where  the behavior of, say, gas molecules in
a certain volume and behavior of groups of active elements (individuals) with
their own interests and sets of actions are considered as the same ones.
The mean field game model is in fact based on the assumption that
a group of human individuals can be replaced
with its generalized representative, whereas all the basic works
in social and economic sciences always emphasized the difference in
behavior of an individual and a group of individuals.

In particular, the well-known K.~Arrow impossibility theorem
 states that a collective preference relation that is
consistent for all the team members is only possible if it matches one of the
individual relations under rather general assumptions
(see, for example, \cite{Eke79,AHS06}). Therefore, the whole team
behavior will be in general more complicated, and it is not represented
by some common preference relationship, i.e. a different model is required instead of
an optimization problem with respect to some preference relation. This assertion
it can be illustrated with simple examples.

Let us consider the election of a team leader when the team consists of $n$
 groups ($n>2$), each $i$-th group nominates its own candidate
$a_{i}$ and sends its representative to the election with its particular
preference relation $\succ_{i}$. The simplest option for this
preference is to consider its own candidate as the best one, i.e. to choose
$$
a_{i}\succ_{i} a_{j}, \ \forall j \neq i.
$$
Then it is impossible to choose the leader, except to simply take one of these
candidates, i.e., to take one of the individual preference relations as the
collective. In this case, the share of its support will be equal to $1/n
\rightarrow 0$ as $n\rightarrow \infty$, i.e. the decision will be absolutely
\lq\lq illegitimate".

It is also well-known that collective  preferences can be non-transitive, even if all
the individual preferences are transitive. It was first noticed by M.~Condorcet in the $XVIII$-th century
(see, for example, \cite{AHS06}). We now give a generalized version of this
paradox. In the previous example, let the representatives compare all the
candidates in pairs as follows:
\begin{eqnarray*}
&& a_{i}\succ_{1} a_{i+1}, \ i=1,\ldots,n-1; \\
&& a_{i}\succ_{k} a_{i+1},  \ i=1,\ldots,k-2,k,\ldots,n-1, \ \mbox{and} \ a_{n}\succ_{k} a_{1},\\
&&  \ \mbox{for} \ k=2,\ldots,n-1;  \\
&& a_{i}\succ_{n} a_{i+1},  \ i=1,\ldots,n-2, \ \mbox{and} \
a_{n}\succ_{n} a_{1}.
\end{eqnarray*}
If we create the collective preference based on the majority
votes rule, we get the cycle
$$
a_{i}\succ a_{i+1}, \ i=1,\ldots,n-1, \ \mbox{and} \ a_{n}\succ a_{1}.
$$
At the same time, the support share for any pairwise preference above
is equal to $(n-1)/n \rightarrow 1$ as $n\rightarrow \infty$.

By themselves, models with a large number of
 active elements, including those obtained as limit ones from game models are fairly well-known (see
\cite{Hil74,ABB90, Mas83,NS83}). But all these models are of equilibrium
type and  not reduced in general to optimization problems. Hence, they require somewhat different
mathematical tools.


\subsection{A long-term model  for management of renewable natural resources}

Since the management of natural resources, including
renewable, is essential for sustainable development of our world,
 the corresponding management tasks are in the spotlight of
many researchers. Note that creation of proper models here requires for
many heterogeneous, but interrelated factors to be taken into account,
in particular, economic, environmental, social ones as well as their outcomes
for a fairly long time period. For this reason, the
corresponding models are often very complex from the mathematical point of view,
and involve uncertain factors. Let us take the
long-term forest management model as an example.

First of all, the forests are used for different purposes. In Economics,
forest is a source of timber and fuel wood, in Ecology,
it is an environment of air purification and carbon absorption,
in Biology, it is the habitat of animals, birds, and plants.
Also, the forest may be a place of rest from the point of view of human society.
Any attempt to combine all these factors into a
single goal (utility) function and to specify all the relationships and restrictions,
including the effects of weather, environmental pollution, and invasions of
harmful insects, etc., will result in an actually unsolvable
task.

It is therefore necessary to make some decomposition of the problem into
several quite independent blocks. For example,  an arrangement of goals
with a decomposition of the planning periods was suggested in \cite{Kon14}.
Since the sustainable development is the principal goal, it seems better to remove
evaluations of economic factors from long-term models.
Note that there are  a lot of models in this area,
which are formulated mainly as optimal control problems, i.e. as  problems
of maximization of quality criteria along the movement trajectories described by a
state equation (see, for example, \cite{AOK09}).

Let us describe this discrete time model; where the time horizon is we divided into stages (years)
$t=1,2,\ldots,T$. The total forest territory
is bounded above by $S$ and is divided into many stands, each stand
containing only trees of the same age $i=1,\ldots,L$. For
simplicity, we suppose that there is only one kind of trees, since
the case of many species only increases dimensionality within the
same model. We denote by $v^{t}=(v^{t}_{1},\ldots,v^{t}_{L})$  the forest territory vector at the
 end of the $t$-th stage, the  initial distribution
$v^{0}$ is supposed to be known.
Next, the natural forest dynamics in the absence of external influences can be described by the difference relation
$$
v^{t+1}= A(v^t), \quad t=0,1,\ldots,T-1;
$$
where the operator $A$ in the simplest case is given by a matrix of the order $L$.
If necessary, it may contain undefined factors, but
this will not affect the essence of the model.
We denote by $u^{t}=(u^{t}_{1},\ldots, u^{t}_{L})$ and $w^{t}=(w^{t}_{1},\ldots,w^{t}_{L})$ the harvest and planting
territory (square) vectors (say, in hectares) within the $t$-th time
stage, respectively, while
$u^{t}_{i} = 0$ for $i=1,\ldots, l-1$, and $w^{t}_{i} = 0$ for
$i=l_{0}+1,\ldots,L$, i.e. $l$ is the minimal harvesting age, $l_{0}$ is the maximal planting age.
Also, $L$ is considered as the maximal product age.
Then the dynamics of the forest will be
described by the relations:
\begin{eqnarray*}
&& v^{t+1}= A (v^t-u^t)+w^t,  \quad t=0,1,\ldots,T-1; \\
&& \sum_{i=1}^{L}v^{t}_{i} \le S,  \
 v^t \geq 0,  \ u^t \geq 0, \ w^t \geq 0,
\quad t=0,1,\ldots,T.
\end{eqnarray*}
At the end of the planning period, some desired set of feasible distributions is given, for example,
$$
v^{T} \in V,
$$
where $V$ is some set in $ \mathbb{R}^{L}$.
Among various additional conditions, we impose only lower bounds
(i.e. minimal feasible volumes) $\Gamma_{t}$ of carbon
sequestration per stage $t$ and insert the constraints
$$
\sum_{i=1}^{L}\gamma_{i}v^{t}_{i}  \ge \Gamma_{t}, \quad
t=1,2,\ldots,T,
$$
where $\gamma_{i}$ denotes the carbon sequestration volume from one hectare of forest of age $i$ per stage, for
$i=1,\ldots,L$.

The set of feasible trajectories $\{ v^t \}$, $ \{ u^t \}$, $\{ w^t \}$,
which satisfy the above restrictions, may be rather large.
Therefore, one should choose the quality criterion for
the management. The traditional approach is to take the profit function
along the trajectory. It suffices to determine the timber price (say, per cubic meter)
and timber yield of each forest age, the unit cost of felling a hectare of forest of
each product age, as well as the unit cost of planting seedlings
 of each suitable age per hectare. Then one
 can calculate the profit value along a given trajectory.
The main drawback of this approach is that the values $l$ and $L$ are
large enough. For the pine forest they are the following: $l \approx 60$ years, $L\approx
120$ years. Then the planning horizon $T$ should be even longer in order to take into account
the full rotation age, hence, $T \approx 200$. Therefore, any price
values for such a long period of time will be unrealistic, and it is necessary to
take a different criterion. For example, one can take a cost comparison based criterion.
Let $\mu_{i}$ denote the yield of timber  from one hectare of forest of age $i$ and $ \eta_{i}$
denote the yield from one hectare of forest of age $i$ after attaining age $l$
(or some other reference age). Then one can take define the goal function
$$
 \sum \limits_{t=1}^{T} \left(\sum
\limits_{i=l}^{L} \mu _{i}u^{t}_{i} -
          \sum \limits_{i=1}^{l_{0}} \eta_{i}w^{t}_{i}\right)
$$
and choose the feasible trajectory that delivers the maximal
 value of this  criterion. If necessary, one can apply uniform
criteria. The resulting optimization problems are solved with
known computational methods.


}
\end{document}